\documentclass[12pt,a4paper]{article}%{amsart}%{article}% %
\usepackage{latexsym,amsmath,amssymb,amsfonts,amstext,amsthm}
\usepackage{url}
\usepackage{graphicx}
\usepackage{color}

\newcounter{num} 
\newcommand{\Fg}[1][]{\thenum}
\begin{document}
\title{The Descartes circles theorem and division by zero calculus
}

\author{ 
Hiroshi Okumura and Saburou Saitoh
}

%\date{\today}

\maketitle

\def\a{\alpha}
\def\b{\beta}

\def\z{\zeta}
\def\rA{r_{\rm A}}

\textbf{Abstract.} 
From the viewpoint of the division by zero $(0/0=1/0=z/0=0)$ and the division 
by zero calculus, we will show that in the very beautiful theorem by Descartes 
on three touching circles is valid for lines and points for circles except for 
one case. However, for the exceptional case, we can obtain interesting results 
from the division by zero calculus.  

\section{Introduction}

We recall the famous and beautiful theorem (\cite{d2002,s1936}):
\medskip

\textsc{\bf Theorem (Descartes)} {\sl Let $C_i$ $(i=1,2,3)$ be 
circles touching to each other of radii $r_i$.  
If a circle $C_4$ touches the three circles, then its radius $r_4$ is 
given by 
\begin{equation}\label{eq1}
\frac{1}{r_4} = \frac{1}{r_1} + \frac{1}{r_2} + \frac{1}{r_3}
 \pm 2 \sqrt{\frac{1}{r_1r_2} + \frac{1}{r_2r_3} + \frac{1}{r_3r_1}}.  
\end{equation}
}

\medskip

As well-known, circles and lines may be looked as  the same ones in complex 
analysis, in the sense of stereographic projection and many reasons.
Therefore, we will consider whether the theorem is valid for line cases and 
point cases for circles. Here, we will discuss this problem clearly from the 
division by zero viewpoint.
The Descartes circle theorem is valid except for  one case for lines and 
points for the three circles and for one exception case, we can obtain very 
interesting results, by the division by zero calculus.

\section{The division by zero calculus}

For any  Laurent expansion around $z=a$,

\begin{equation}\label{eq2}
f(z) = \sum_{n=-\infty}^{\infty} C_n (z - a)^n,  
\end{equation}
we obtain the identity, by the division by zero

\begin{equation}\label{eq3}
f(a) =  C_0. 
\end{equation}
(Here, as convention, we consider as $0^0=1$.)
\medskip

For the correspondence \eqref{eq3} for the function $f(z)$, we will call it 
{\bf the division by zero calculus}. By considering the derivatives in 
\eqref{eq2}, we can define any order derivatives of the function $f$ at the 
singular point $a$.

We have considered our mathematics around an isolated singular point for 
analytic functions, however, we do not consider mathematics at the singular 
point itself. At the isolated singular point, we consider our mathematics with 
the limiting concept, however, the limiting values to the singular point and 
the value at the singular point of the function are different.
By the division by zero calculus, we can consider the values and differential 
coefficients at the singular point. 

The division by zero ($0/0=1/0=z/0=0$) is trivial and clear in the natural 
sense of the generalized division (fraction) against its mysterious and long 
history (see for example, \cite{romig}), since we know the Moore-Penrose 
generalized inverse for the elementary equation $az =b$. Therefore, the 
division by zero calculus above and its applications are important. See the 
references \cite{s14,kmsy,msy15,s16,ms16,ms18,mos17,osm17, mms18,ps18} for 
the details and the related topics.
We regret that our common sense for the division by zero are still wrong; one 
typical comment for our division by zero results is given by some physician:

{\it Here is how I see the problem with prohibition on division by zero, 
which is the biggest scandal in modern mathematics as you rightly pointed 
out} (2017.10.14.8:55).

\medskip

However, in this paper we do not need any information and results in the 
division by zero, we need only the definition \eqref{eq3} of the division by 
zero calculus.

As stated already in \cite{mos17}, in general, for a circle with radius $r$, 
its curvature is given by $1/r$ and by the division by zero, for the point 
circle, its curvature is zero. Meanwhile, for a line corresponding the case 
$r = \infty$, its curvature is also zero, however, then we should consider the 
case as $r = 0$, not $\infty$. For this reality and reasonable situation, look 
the paper. By this interpretation, we will show that the theorem is valid for 
line and point circle cases clearly. We would like to show that the division 
by zero  ($0/0=1/0=z/0=0$) is very natural also from the viewpoint of the 
Descartes circles theorem.

\section{Results}

We would like to consider all the cases for the Descartes theorem for lines 
and point circles, step by step.

\subsection{One line and two circles case}

We consider the case in which the circle $C_3$ is one of the external common  
tangents of the circles $C_1$ and $C_2$. This is a typical case in this paper. 
We assume $r_1\ge r_2$. 
We now have $r_3=0$ in \eqref{eq1}. Hence 
$$
\frac{1}{r_4} = \frac{1}{r_1} + \frac{1}{r_2} + \frac{1}{0}
\pm 2\sqrt{\frac{1}{r_1r_2} + \frac{1}{r_2\cdot 0} + \frac{1}{0\cdot r_1}}
= \frac{1}{r_1} + \frac{1}{r_2} \pm 2 \sqrt{\frac{1}{r_1r_2}}. 
$$
This implies  
$$
\frac{1}{\sqrt{r_4}}=\frac{1}{\sqrt{r_1}}+\frac{1}{\sqrt{r_2}}
$$
in the plus sign case. The circle $C_4$ is the incircle of the curvilinear 
triangle made by $C_1$, $C_2$ and $C_3$  (see Figure \ref{fin}). 
In the minus sign case we have 
$$
\frac{1}{\sqrt{r_4}}=\frac{1}{\sqrt{r_2}}-\frac{1}{\sqrt{r_1}}. 
$$
In this case $C_2$ is the incircle of the curvilinear 
triangle made by the other three (see Figure \ref{fex}).

\medskip
\begin{minipage}{.42\hsize}
\begin{center}
\includegraphics[clip,width=67mm]{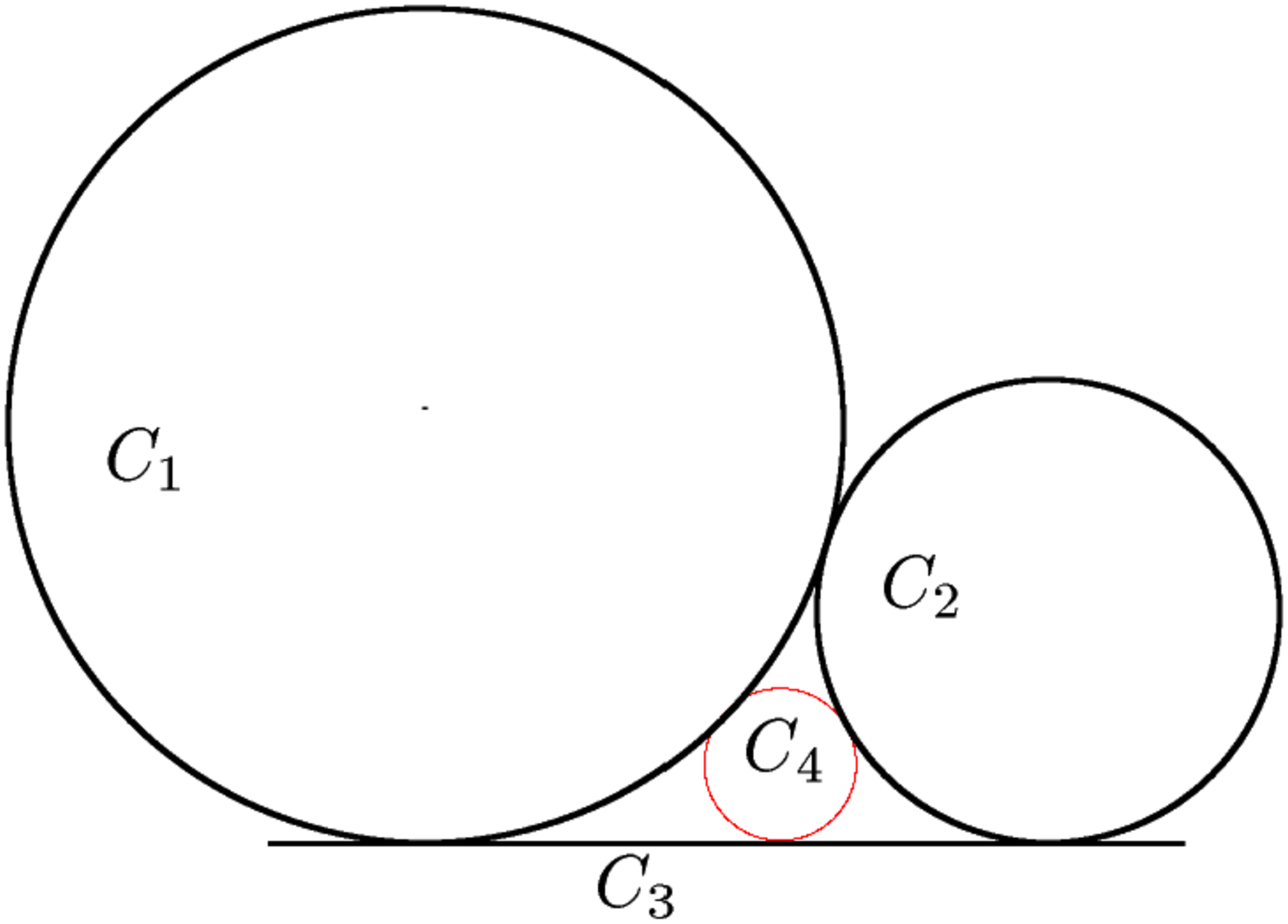}\refstepcounter{num}\label{fin}\\
Figure \Fg . 
\end{center}
\end{minipage}
\begin{minipage}{.52\hsize}
\begin{center}
\includegraphics[clip,width=62mm]{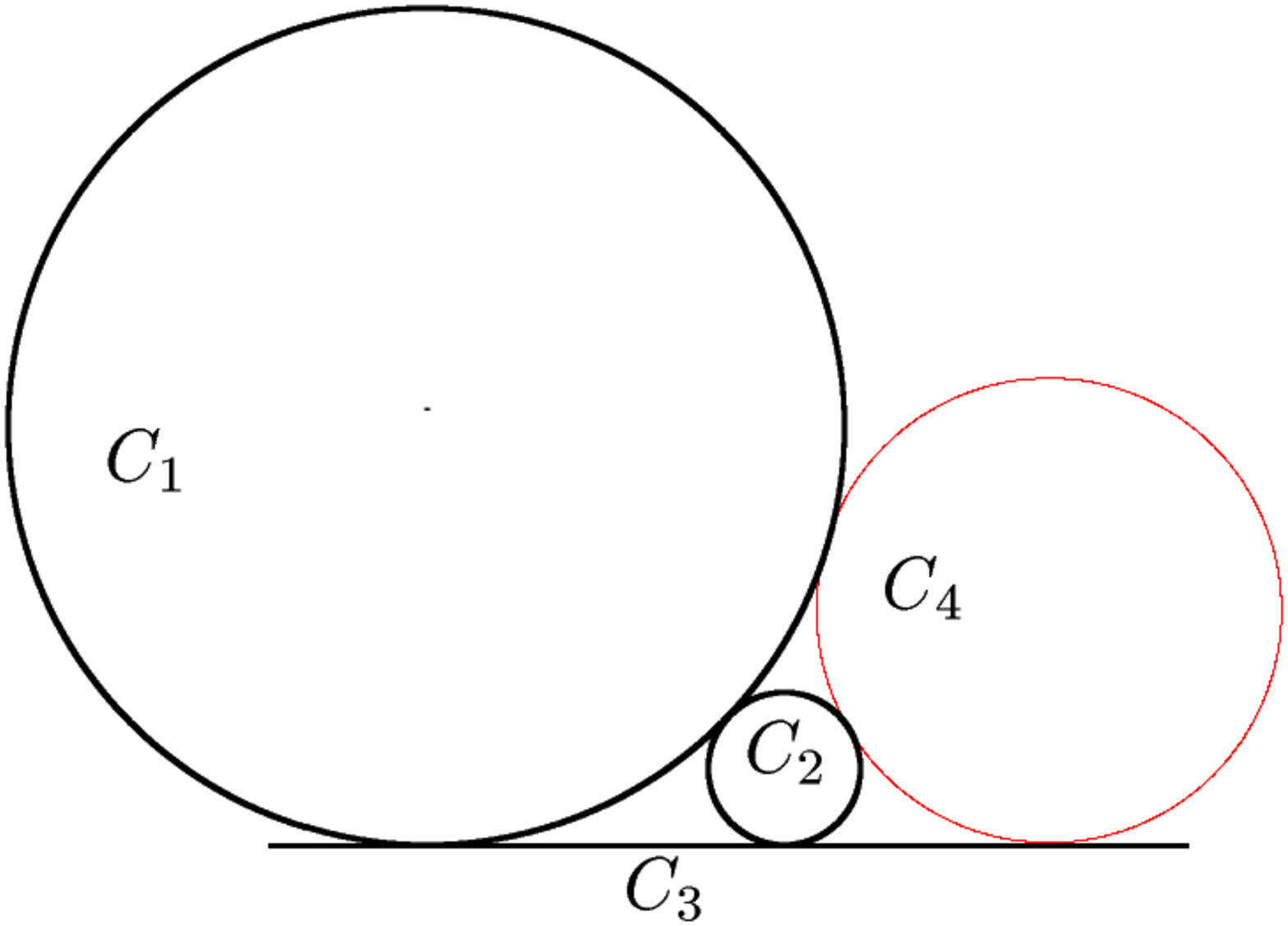}\refstepcounter{num}\label{fex}\\
Figure \Fg . 
\end{center}
\end{minipage}

\medskip
Of course, the result is known. The result was also well-known in 
Wasan geometry \cite{ucd} with the Descartes circle theorem itself 
\cite{zs}, \cite{6}.

\subsection{Two lines and one circle case}

In this case, the two lines have to be parallel, and so, this case is trivial, 
because then the two circles are congruent, by the division by 
zero $1/0=0$.

\subsection{One point circle and two circles case}

This case is another typical case for the theorem. Intuitively, for $r_3 =0$,  the circle $C_3$ is  the common point of the circles 
$C_1$ and $C_2$. Then, there does not exist any touching circle of the three 
circles $C_j; j =1,2,3$.

For the point circle $C_3$, we will consider it by limiting of circles 
touching the circles $C_1$ and $C_2$ to the common point. Then, we will examine 
the circles $C_4$ and the Descartes theorem.

Let us consider the circles $C_1 : (x + r_1)^2 +y^2 = r_1^2$ and 
$C_2: (x - r_2)^2 + y^2 = r_2^2$ ($r_i>0)$.  
We recall the parametric representation of any such circle $C_3$: H. Okumura 
and M. Watanabe gave a theorem in \cite{ow2004}, for a real number $z$, the 
point $(0, 2\sqrt{r_1 r_2}/z)$ is denoted by $V_z$, then 
\medskip

\textsc{\bf Theorem 7.} {\sl The circle $C_3$ passing through $V_{z\pm1}$ for 
a real number 
$z\not=\pm1$ and touching the circles $C_1$ and $C_2$ can be represented 
by the equation 

$$
\left(x-{r_1-r_2\over z^2-1}\right)^2+
\left(y-{2z\sqrt{r_1 r_2}\over z^2-1}\right)^2
=\left({r_1+r_2\over z^2-1}\right)^2. 
$$
}
\medskip

By setting $z = 1/w$, we will consider the case $w=0$; that is, the case 
$z=\infty$ in the classical sense; that is, the circle $C_3$ is reduced to 
the origin.

We look for the circle $C_4$ touching the three circles $C_j; j= 1,2,3$. 
We set
\begin{equation}\label{eqcircle} %old(5)
C_4: (x - x_4)^2 + (y - y_4)^2 = r_4^2.
\end{equation}

Then, from the touching property we obtain:
$$
x_4=\frac{r_1 r_2(r_1 - r_2)w^2}{D},
$$
$$
y_4=\frac{2r_1 r_2\left(\sqrt{r_1 r_2}+(r_1 + r_2)w\right)w}{D} 
$$
and
$$
r_4=\frac{r_1 r_2(r_1 + r_2)w^2}{D}, 
$$
where
$$
D=r_1 r_2+2\sqrt{r_1 r_2}(r_1 + r_2)w+(r_1^2+ r_1 r_2 +r_2^2)w^2.
$$
Notice that there are four sets of the solutions of $x_4$, $y_4$, 
$r_4$, but we consider only one set, because the other cases can be 
considered similarly. 

By inserting these values to \eqref{eqcircle}, we obtain
$$
f_0+f_1w+f_2w^2=0, 
$$
where
$$
f_0=r_1 r_2(x^2+y^2),
$$
$$
f_1=2\sqrt{r_1 r_2}\big((r_1 + r_2)(x^2+y^2)-2r_1 r_2y\big),
$$
and
$$
f_2=(r_1^2+ r_1 r_2 +r_2^2)(x^2+y^2)+2r_1 r_2(r_2-r_1)x-4r_1r_2(r_1 + r_2)y+4r_1^2r_2^2.
$$
By using the division by zero calculus for $w = 0$, we obtain, for the first, 
for $w=0$, the second by setting $w=0$ after dividing by  $w$  and for the 
third case, by setting $w=0$ after dividing by $w^2$, 

\begin{equation}\label{poincircle}%old(7)
x^2+y^2=0   
\end{equation}
\begin{equation}\label{bankoff}%old(8)
(r_1+ r_2)(x^2+y^2)-2r_1 r_2y=0 
\end{equation}
\begin{equation}\label{incircle}%old(9)
(r_1^2+r_1 r_2+r_2^2)(x^2+y^2)+2r_1 r_2(r_2 -r_1)x-4r_1 r_2(r_1 +  r_2)y+4r_1^2r_2^2=0. 
\end{equation}

Note that \eqref{bankoff} is the red circle in Figure \ref{farb} and its 
radius is
\begin{equation}\label{redcirc1}%old(10)
\dfrac{r_1 r_2}{r_1 +r_2}
\end{equation}
and \eqref{incircle} is the green circle in Figure \ref{farb} whose radius 
is 
\begin{equation}\label{inrad}%old(11)
\dfrac{r_1 r_2(r_1 +r_2)}{r_1^2+ r_1 r_2+r_2^2}. 
\end{equation}

\medskip
\begin{center}
\includegraphics[clip,width=100mm]{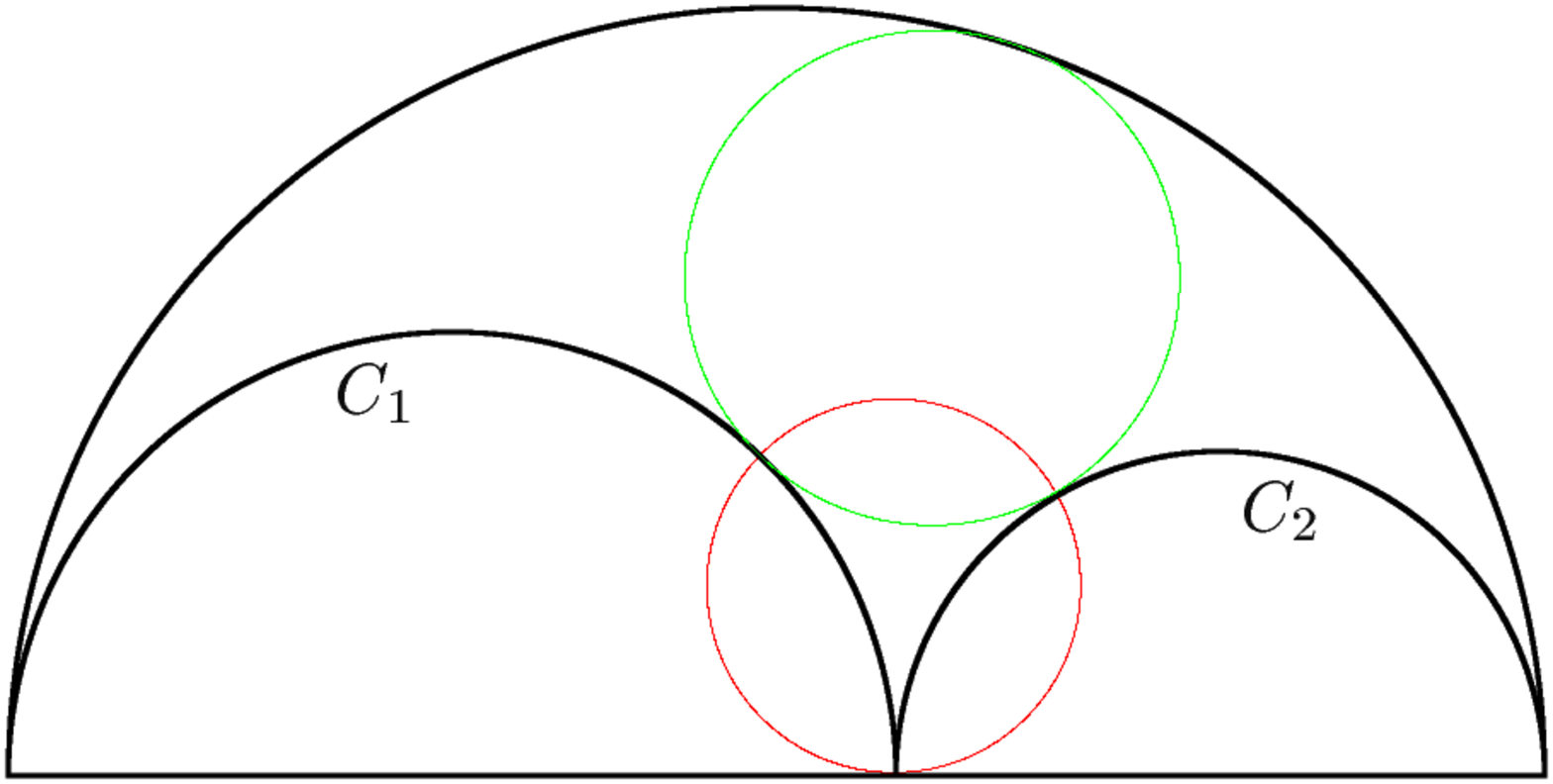}\refstepcounter{num}\label{farb}\\
Figure \Fg . 
\end{center}

When the circle $C_3$ is reduced to the origin, of course, the inscribed 
circle $C_4$ is reduced to the origin, then the Descartes theorem is not valid.
However, by the division by zero calculus, then the origin of $C_4$ is changed 
suddenly for the cases \eqref{poincircle}, \eqref{bankoff} and 
\eqref{incircle}, and for the circle \eqref{bankoff}, the Descartes theorem is 
valid for $r_3=0$, surprisingly.

Indeed, in \eqref{eq1} we set $\xi = \sqrt{r_3}$, then \eqref{eq1} is as follows:
$$
\frac{1}{r_4} = \frac{1}{r_1} + \frac{1}{r_2} + \frac{1}{\xi^2}
\pm 2  \frac{1}{\xi} \sqrt{\frac{\xi^2}{r_1r_2} + \left(\frac{1}{r_1} 
+ \frac{1}{r_2}\right)}.
$$
and so, by the division by zero calculus at $\xi =0$, we have
$$
\frac{1}{r_4} = \frac{1}{r_1} + \frac{1}{r_2}
$$
which is  \eqref{redcirc1}. Note, in particular, that the division by zero 
calculus may be applied in many ways and so, for the results obtained should 
be examined some meanings.
This circle \eqref{bankoff} may be looked a circle touching the origin and 
two circles 
$C_1$ and $C_2$, because by the division by zero calculus
$$
\tan \frac{\pi}{2} = 0,
$$
that is a popular property (see \cite{mos17}).

%Let $C_0$ be the circle touching $C_1$ and $C_2$ 
%at the points on the 
%line joining the centers of $C_1$ and $C_2$

Meanwhile, the circle \eqref{incircle} touches the 
circles $C_1$,  $C_2$ and the beautiful circle  with center $(r_2-r_1,0)$ 
with radius $r_1 + r_2$. The each of the areas surrounded by the last circle 
and the circles $C_1$, $C_2$ is called an arbelos, and 
the circle \eqref{bankoff} is the famous Bankoff circle of the arbelos. 

For $r_3=-(r_1 +r_2)$, from the Descartes identity \eqref{eq1}, $r_4$ equals 
\eqref{inrad}. That 
is, when we consider that the circle $C_3$ is changed to the circle with 
center $(r_2-r_1,0)$ with radius $r_1 + r_2$, the Descartes identity 
holds. Here, the minus sign shows that the circles $C_1$ and $C_2$ touch $C_3$ 
internally from the inside of $C_3$. 

\subsection{Two point circles and one circle case}

This case is trivial, because the externally touching circle 
is coincident with one circle.

\subsection{Three points case and three lines case}

In these cases we have  $r_j = 0,  j=1,2,3$ and the formula \eqref{eq1} shows 
that $r_4= 0$. This statement is trivial in the general sense.

As the solution of the simplest equation
\begin{equation}\label{ax=b}%old(12)
ax =b,
\end{equation}
we have $x=0$ for $a=0, b\ne 0$ as the standard value, or the Moore-Penrose 
generalized inverse. This will mean in a sense, the solution does not exist; 
to solve the equation \eqref{ax=b} is impossible.
We saw for different parallel lines or different parallel planes, their common 
points are the origin in \cite{mos17}. Certainly they have the common point of 
the point at infinity and the point at infinity is represented by zero. 
However, we can understand also that they have  no common points, because the 
point at infinity is an ideal point. The zero will represent some impossibility.

In the Descartes theorem, three lines and three points cases, we can 
understand that the touching circle does not exist, or it is the point and 
so  the Descartes theorem is valid.

\section{Conclusion}
By the division by zero calculus, we were able to give a general theorem of 
Descartes for containing lines and point circles for circles, simply. At the 
same time, we showed that the division by zero calculus is natural for the 
famous Descartes circles theorem clearly.

\bibliographystyle{plain}

\begin{thebibliography}{99} 

\bibitem{d2002}
C. Jeffrey, C.L. Lagarias,  A. R. Mallows and A. R. Wilks,  
Beyond the Descartes Circle Theorem. The American Mathematical Monthly 
{\bf 109}(4) (2002) 338--361. doi:10.2307/2695498. JSTOR 2695498.

\bibitem{kmsy}
M. Kuroda, H. Michiwaki, S. Saitoh and M. Yamane,
New meanings of the division by zero and interpretations on $100/0=0$ and 
on $0/0=0$, 
Int. J. Appl. Math. {\bf 27}(2) (2014), 191-198,  DOI: 10.12732/ijam.v27i2.9.

\bibitem{ms16}
T. Matsuura and S. Saitoh, Matrices and division by zero�$z/0=0$, 
Advances in Linear Algebra \& Matrix Theory, \textbf{6}(2) 2016, 51--58
Published Online June 2016 in SciRes. 
\url{http://www.scirp.org/journal/alamt}. 
\url{http://dx.doi.org/10.4236/alamt.2016.62007}.

\bibitem{ms18}
T. Matsuura and S. Saitoh,
Division by zero calculus and singular integrals. 
(Submitted for publication).

\bibitem{mms18}
T. Matsuura, H. Michiwaki and S. Saitoh,
$\log 0= \log \infty =0$ and applications. Differential and Difference 
Equations with Applications. Springer Proceedings in Mathematics \& Statistics.

\bibitem{msy15}
H. Michiwaki, S. Saitoh and  M.Yamada, 
Reality of the division by zero $z/0=0$. 
IJAPM  International J. of Applied Physics and Math. \textbf{6} (2015), 1--8. 
\url{http://www.ijapm.org/show-63-504-1.html}. 

\bibitem{mos17}
H. Michiwaki, H. Okumura and S. Saitoh,
 Division by Zero $z/0 = 0$ in Euclidean Spaces,
 International Journal of Mathematics and Computation, \textbf{28}(1) (2017) 
 1--16.
 
\bibitem{ow2004}
H. Okumura and T. Watanabe,
The Twin Circles of Archimedes in a Skewed Arbelos,
Forum Geometricorum,  \textbf{4} (2004), 229--251.  

\bibitem{osm17}
H. Okumura, S. Saitoh and T. Matsuura, Relations of $0$ and $\infty$,
Journal of Technology and Social Science (JTSS), \textbf{1} (2017),  70--77.

\bibitem{ps18}
S. Pinelas and S. Saitoh,
Division by zero calculus and differential equations, Differential and 
Difference Equations with Applications. 
Springer Proceedings in Mathematics \& Statistics.

\bibitem{romig}
H. G. Romig, Discussions: Early History of Division by Zero,
American Mathematical Monthly, \textbf{31}(8) (Oct., 1924), 387--389.

\bibitem{s14}
S. Saitoh, Generalized inversions of Hadamard and tensor products for matrices,  Advances in Linear Algebra \& Matrix Theory. {\bf 4}(2) (2014) 87--95.
\url{http://www.scirp.org/journal/ALAMT/}. 

\bibitem{s16}
S. Saitoh, A reproducing kernel theory with some general applications, 
Qian,T./Rodino,L.(eds.): Mathematical Analysis, Probability and 
Applications - Plenary Lectures: Isaac 2015, Macau, China, Springer 
Proceedings in Mathematics and Statistics,  {\bf 177} (2016), 151-182. 
(Springer). 

\bibitem{s1936}
 F. Soddy, The Kiss Precise. Nature \textbf{137}(3477) (1936), 1021. 
doi:10.1038/1371021a0.

\bibitem{6} Sugano or Kanno, Rokusha Jutsu, 1798, 
Tohoku University Wasan Material Database, 
\url{http://www.i-repository.net/il/meta_pub/G0000398wasan_4100007297}. 

\bibitem{ucd} Uchida ed., Z\=oho Sany\=o Tebikigusa, 1764, 
Tohoku University Wasan Material Database, 
\url{http://www.i-repository.net/il/meta_pub/G0000398wasan_4100005700}.

\bibitem{zs} Yamaji, Zeishiki Endan, 1751, 
Tohoku University Wasan Material Database,
\url{http://www.i-repository.net/il/meta_pub/G0000398wasan_4100001519}. 

\end{thebibliography}

\noindent
Hiroshi Okumura,\\ 
Maebashi Gunma 371-0123, Japan\\
{\it  E-mail address}: okmr@protonmail.com
\medskip

\noindent
Saburou Saitoh, \\ 
Institute of Reproducing Kernels\\
Kawauchi-cho, 5-1648-16, Kiryu 376-0041, Japan\\
{\it  E-mail address}: kbdmm360@yahoo.com.jp

\end{document}